\newif\iftopology
\topologyfalse
\iftopology
\documentclass{elsart}
\else
\documentclass[twocolumn, 11pt]{article}
\usepackage{mathptmx}
\usepackage[colorlinks=true, linkcolor=blue, citecolor=blue, urlcolor=blue]{hyperref}
\fi

\newcommand{\error}{\textsuperscript{*}}
\newcommand{\errortwoalone}{\textsuperscript{$\dagger$}}
\newcommand{\errortwo}{\textsuperscript{*$\dagger$}}

\makeatletter
\newenvironment{semitable}
               {\small
                 \let\\\@centercr
                 \list{}{\itemsep      \z@
                        \itemindent   -1.5em%
                        \listparindent\itemindent
                        \rightmargin  0.0em
                        }%
                \item\relax}
               {\endlist}
\makeatother

\begin{document}

\iftopology
\else

    \begin{center}
   {\Large \textbf{A table of boundary slopes of \\
    Montesinos knots\\}}

    \bigskip

    {\large  Nathan M. Dunfield\footnote{Partially supported by
    an NSF Graduate Fellowship and a Sloan Dissertation Fellowship.}}

     \medskip

     {Revised version of June 2010}

     \medskip

     {email: \sl nathan@dunfield.info}

     \end{center}

    \begin{abstract}
    This note corrects errors in Hatcher and Oertel's table of boundary 
     slopes of Montesinos knots which have projections with 10 or fewer
     crossings.
     \end{abstract}

\fi

\section{Introduction}

In \cite{HatcherOertel}, Hatcher and Oertel gave an algorithm for
computing the boundary slopes of a Montesinos knot.  At the end of
their paper they provided a table giving the boundary slopes for each
Montesinos knot in the standard table in Rolfsen's book
\cite{Rolfsen76}.  Unfortunately, their table contains several (17)
errors.  These errors were due to problems with the computer program
that generated the table, as well as transcription/printing errors.
This note provides a corrected table which was generated by a
completely new computer program.  Section~\ref{table} contains the
corrected table.  Section~\ref{method} describes the precautions I
took in writing the new program.  Unfortunately, not all mistakes in
\cite{HatcherOertel} were found in the original
\href{http://arxiv.org/abs/math.GT/9901120v2}{1999 version of this
  note}, and this 2010 revision corrects several more (details are
given in Section~\ref{method}).   Section~\ref{gettingprogram} describes
where the reader can download the new program.

I thank Allen Hatcher and especially Ulrich Oertel for their help with
this project.  I also thank Thomas Mattman and Marc Culler their help
with finding and correcting the additional errors.

\section{Revised table}\label{table}

Throughout this section \error\ denotes an entry which differs from
the corresponding entry in the table in \cite{HatcherOertel}, and
\errortwoalone\ denotes one that differs from the 1999 version of this note. 

Below are the boundary slopes for those Montesinos knots of ten or
fewer crossings.  Notation follows \cite{HatcherOertel} and
\cite{Rolfsen76}.

\begin{semitable}
$3_{1} = \mathrm{K}(1/3):$  $0$, $6$\\ 
$4_{1} = \mathrm{K}(2/5):$  $-4$, $0$, $4$\\ 
$5_{1} = \mathrm{K}(1/5):$  $0$, $10$\\ 
$5_{2} = \mathrm{K}(3/7):$  $0$, $4$, $10$\\ 
$6_{1} = \mathrm{K}(4/9):$  $-4$, $0$, $8$\\ 
$6_{2} = \mathrm{K}(4/11):$  $-4$, $0$, $2$, $8$\\ 
$6_{3} = \mathrm{K}(5/13):$  $-6$, $-2$, $0$, $2$, $6$\\ 
$7_{1} = \mathrm{K}(1/7):$  $0$, $14$\\ 
$7_{2} = \mathrm{K}(5/11):$  $0$, $4$, $14$\\ 
$7_{3} = \mathrm{K}(4/13):$  $-14$, $-8$, $0$\\ 
$7_{4} = \mathrm{K}(4/15):$  $-14$, $-8$, $0$\\ 
$7_{5} = \mathrm{K}(7/17):$  $0$, $4$, $6$, $10$, $14$\\ 
$7_{6} = \mathrm{K}(7/19):$  $-4$, $0$, $4$, $6$, $10$\\ 
$7_{7} = \mathrm{K}(8/21):$  $-8$, $-4$, $0$, $6$\\ 
$8_{1} = \mathrm{K}(6/13):$  $-4$, $0$, $12$\\ 
$8_{2} = \mathrm{K}(6/17):$  $-4$, $0$, $6$, $12$\\ 
$8_{3} = \mathrm{K}(4/17):$  $-8$, $0$, $8$\\ 
$8_{4} = \mathrm{K}(5/19):$  $-8$, $-2$, $0$, $8$\\ 
$8_{5} = \mathrm{K}(1/3, 1/3, 1/2):$  $-4$, $0$, $2$, $8$, $10$, $12$\\ 
$8_{6} = \mathrm{K}(10/23):$  $-4$, $0$, $2$, $6$, $12$\\ 
$8_{7} = \mathrm{K}(9/23):$  $-10$, $-6$, $-2$, $0$, $6$\\ 
$8_{8} = \mathrm{K}(9/25):$  $-10$, $-6$, $-4$, $0$, $2$, $6$\\ 
$8_{9} = \mathrm{K}(7/25):$  $-8$, $-2$, $0$, $2$, $8$\\ 
$8_{10} = \mathrm{K}(1/3, 2/3, 1/2):$  $-6$, $-2$, $0$, $6$, $8$, $10$\\ 
$8_{11} = \mathrm{K}(10/27):$  $-4$, $0$, $6$, $12$\\ 
$8_{12} = \mathrm{K}(12/29):$  $-8$, $-4$, $0$, $4$, $8$\\ 
$8_{13} = \mathrm{K}(11/29):$  $-10$, $-6$, $-4$, $-2$, $0$, $6$\\ 
$8_{14}\error = \mathrm{K}(12/31):$  $-4$, $0$, $4$, $6$, $8$, $12$\\ 
$8_{15} = \mathrm{K}(2/3, 2/3, 1/2):$  $-16$, $-12$, $-10$, $-8$, $-4$, $-2$, $0$\\ 
$8_{19} = \mathrm{K}(1/3, 1/3, -1/2):$  $0$, $12$\\ 
$8_{20} = \mathrm{K}(1/3, 2/3, -1/2):$  $-10$, $0$, $8/3$\\ 
$8_{21} = \mathrm{K}(2/3, 2/3, -1/2):$  $-12$, $-6$, $-2$, $0$, $1$\\ 
$9_{1} = \mathrm{K}(1/9):$  $0$, $18$\\ 
$9_{2} = \mathrm{K}(7/15):$  $0$, $4$, $18$\\ 
$9_{3} = \mathrm{K}(6/19):$  $-18$, $-12$, $0$\\ 
$9_{4} = \mathrm{K}(5/21):$  $0$, $8$, $18$\\ 
$9_{5} = \mathrm{K}(6/23):$  $-18$, $-12$, $-8$, $0$\\ 
$9_{6} = \mathrm{K}(11/27):$  $0$, $4$, $10$, $14$, $18$\\ 
$9_{7} = \mathrm{K}(13/29):$  $0$, $4$, $6$, $10$, $18$\\ 
$9_{8} = \mathrm{K}(11/31):$  $-8$, $-4$, $0$, $4$, $6$, $10$\\ 
$9_{9} = \mathrm{K}(9/31):$  $0$, $6$, $8$, $14$, $18$\\ 
$9_{10} = \mathrm{K}(10/33):$  $-18$, $-12$, $-6$, $0$\\ 
$9_{11} = \mathrm{K}(14/33):$  $-14$, $-10$, $-4$, $0$, $4$\\ 
$9_{12} = \mathrm{K}(13/35):$  $-4$, $0$, $6$, $8$, $14$\\ 
$9_{13} = \mathrm{K}(10/37):$  $-18$, $-14$, $-12$, $-8$, $-6$, $0$\\ 
$9_{14} = \mathrm{K}(14/37):$  $-12$, $-8$, $-4$, $0$, $6$\\ 
$9_{15} = \mathrm{K}(16/39):$  $-14$, $-10$, $-8$, $-4$, $0$, $4$\\ 
$9_{16} = \mathrm{K}(1/3, 1/3, 3/2):$  $0$, $4$, $6$, $10$, $12$, $14$, $16$, $18$\\ 
$9_{17} = \mathrm{K}(14/39):$  $-8$, $-4$, $-2$, $0$, $4$, $10$\\ 
$9_{18} = \mathrm{K}(17/41):$  $0$, $4$, $8$, $10$, $12$, $14$, $18$\\ 
$9_{19} = \mathrm{K}(16/41):$  $-8$, $-4$, $0$, $4$, $10$\\ 
$9_{20} = \mathrm{K}(15/41):$  $-4$, $0$, $2$, $6$, $8$, $14$\\ 
$9_{21} = \mathrm{K}(18/43):$  $-14$, $-10$, $-8$, $-4$, $0$, $4$\\ 
$9_{22} = \mathrm{K}(3/5, 1/3, 1/2):$  $-8$, $-4$, $-2$, $0$, $2$, $4$, $6$, $8$, $10$\\ 
$9_{23} = \mathrm{K}(19/45):$  $0$, $4$, $8$, $10$, $14$, $18$\\ 
$9_{24} = \mathrm{K}(1/3, 2/3, 3/2):$  $-10$, $-6$, $-4$, $0$, $2$, $4$, $6$, $8$\\ 
$9_{25} = \mathrm{K}(2/5, 2/3, 1/2):$  $-14$, $-10$, $-8$, $-6$, $-4$, $-2$, $0$, $2$, $4$\\ 
$9_{26} = \mathrm{K}(18/47):$  $-12$, $-8$, $-6$, $-4$, $0$, $6$\\ 
$9_{27} = \mathrm{K}(19/49):$  $-8$, $-4$, $-2$, $0$, $2$, $4$, $6$, $10$\\ 
$9_{28} = \mathrm{K}(2/3, 2/3, 3/2):$  $-12$, $-8$, $-6$, $-2$, $0$, $2$, $4$, $6$\\ 
$9_{30} = \mathrm{K}(3/5, 2/3, 1/2):$  $-10$, $-6$, $-4$, $-2$, $0$, $2$, $4$, $6$, $8$\\ 
$9_{31} = \mathrm{K}(21/55):$  $-6$, $-2$, $0$, $2$, $6$, $12$\\ 
$9_{35} = \mathrm{K}(1/3, 1/3, 1/3):$  $-18$, $-12$, $-4$, $0$\\ 
$9_{36} = \mathrm{K}(2/5, 1/3, 1/2):$  $-4$, $0$, $2$, $4$, $6$, $8$, $10$, $12$, $14$\\ 
$9_{37} = \mathrm{K}(1/3, 2/3, 2/3):$  $-10$, $-4$, $0$, $4$, $8$\\ 
$9_{42} = \mathrm{K}(2/5, 1/3, -1/2):$  $-8$, $0$, $8/3$, $6$\\ 
$9_{43} = \mathrm{K}(3/5, 1/3, -1/2):$  $-4$, $0$, $6$, $8$, $32/3$\\ 
$9_{44} = \mathrm{K}(2/5, 2/3, -1/2):$  $-10$, $-2$, $0$, $1$, $2$, $14/3$\\ 
$9_{45} = \mathrm{K}(3/5, 2/3, -1/2):$  $-14$, $-10$, $-8$, $-4$, $-2$, $0$, $1$\\ 
$9_{46} = \mathrm{K}(1/3, 1/3, -1/3):$  $-12$, $0$, $2$\\ 
$9_{48} = \mathrm{K}(2/3, 2/3, -1/3):$  $-4$, $0$, $4$, $8$, $11$\\ 
$10_{1} = \mathrm{K}(8/17):$  $-4$, $0$, $16$\\ 
$10_{2} = \mathrm{K}(8/23):$  $-4$, $0$, $10$, $16$\\ 
$10_{3} = \mathrm{K}(6/25):$  $-8$, $0$, $12$\\ 
$10_{4} = \mathrm{K}(7/27):$  $-12$, $-6$, $0$, $8$\\ 
$10_{5} = \mathrm{K}(13/33):$  $-14$, $-10$, $-6$, $0$, $6$\\ 
$10_{6} = \mathrm{K}(16/37):$  $-4$, $0$, $6$, $10$, $16$\\ 
$10_{7} = \mathrm{K}(16/43):$  $-4$, $0$, $6$, $10$, $16$\\ 
$10_{8} = \mathrm{K}(6/29):$  $-8$, $0$, $2$, $12$\\ 
$10_{9} = \mathrm{K}(11/39):$  $-12$, $-6$, $-2$, $0$, $8$\\ 
$10_{10} = \mathrm{K}(17/45):$  $-14$, $-10$, $-6$, $-4$, $0$, $6$\\ 
$10_{11} = \mathrm{K}(13/43):$  $-8$, $-2$, $0$, $6$, $12$\\ 
$10_{12} = \mathrm{K}(17/47):$  $-14$, $-10$, $-8$, $-2$, $0$, $6$\\ 
$10_{13} = \mathrm{K}(22/53):$  $-8$, $-4$, $0$, $4$, $8$, $12$\\ 
$10_{14} = \mathrm{K}(22/57):$  $-4$, $0$, $4$, $8$, $10$, $12$, $16$\\ 
$10_{15} = \mathrm{K}(19/43):$  $-10$, $-6$, $0$, $2$, $4$, $10$\\ 
$10_{16} = \mathrm{K}(14/47):$  $-12$, $-6$, $-2$, $0$, $4$, $8$\\ 
$10_{17} = \mathrm{K}(9/41):$  $-10$, $-2$, $0$, $2$, $10$\\ 
$10_{18} = \mathrm{K}(23/55):$  $-8$, $-4$, $0$, $2$, $4$, $8$, $12$\\ 
$10_{19} = \mathrm{K}(14/51):$  $-10$, $-4$, $-2$, $0$, $2$, $4$, $10$\\ 
$10_{20} = \mathrm{K}(16/35):$  $-4$, $0$, $2$, $6$, $16$\\ 
$10_{21} = \mathrm{K}(16/45):$  $-4$, $0$, $4$, $10$, $16$\\ 
$10_{22} = \mathrm{K}(13/49):$  $-12$, $-6$, $0$, $2$, $8$\\ 
$10_{23} = \mathrm{K}(23/59):$  $-14$, $-10$, $-8$, $-6$, $-4$, $0$, $6$\\ 
$10_{24} = \mathrm{K}(24/55):$  $-4$, $0$, $4$, $6$, $10$, $16$\\ 
$10_{25} = \mathrm{K}(24/65):$  $-4$, $0$, $2$, $6$, $10$, $12$, $16$\\ 
$10_{26} = \mathrm{K}(17/61):$  $-12$, $-6$, $-2$, $0$, $4$, $8$\\ 
$10_{27} = \mathrm{K}(27/71):$  $-14$, $-10$, $-6$, $-4$, $0$, $6$\\ 
$10_{28} = \mathrm{K}(19/53):$  $-14$, $-10$, $-8$, $-4$, $-2$, $0$, $6$\\ 
$10_{29} = \mathrm{K}(26/63):$  $-8$, $-4$, $-2$, $0$, $2$, $4$, $8$, $12$\\ 
$10_{30} = \mathrm{K}(26/67):$  $-4$, $0$, $4$, $8$, $10$, $12$, $16$\\ 
$10_{31} = \mathrm{K}(25/57):$  $-10$, $-6$, $-4$, $0$, $2$, $4$, $10$\\ 
$10_{32} = \mathrm{K}(29/69):$  $-8$, $-4$, $-2$, $0$, $2$, $4$, $6$, $8$, $12$\\ 
$10_{33} = \mathrm{K}(18/65):$  $-10$, $-4$, $-2$, $0$, $2$, $4$, $10$\\ 
$10_{34} = \mathrm{K}(13/37):$  $-14$, $-10$, $-4$, $0$, $2$, $6$\\ 
$10_{35} = \mathrm{K}(20/49):$  $-12$, $-8$, $-4$, $0$, $4$, $8$\\ 
$10_{36} = \mathrm{K}(20/51):$  $-4$, $0$, $4$, $8$, $10$, $16$\\ 
$10_{37} = \mathrm{K}(23/53):$  $-10$, $-6$, $-4$, $0$, $4$, $6$, $10$\\ 
$10_{38} = \mathrm{K}(25/59):$  $-4$, $0$, $4$, $6$, $8$, $10$, $12$, $16$\\ 
$10_{39} = \mathrm{K}(22/61):$  $-4$, $0$, $2$, $4$, $6$, $8$, $10$, $12$, $16$\\ 
$10_{40} = \mathrm{K}(29/75):$  $-14$, $-10$, $-6$, $-4$, $-2$, $0$, $2$, $6$\\ 
$10_{41} = \mathrm{K}(26/71):$  $-8$, $-4$, $0$, $2$, $4$, $6$, $8$, $12$\\ 
$10_{42} = \mathrm{K}(31/81):$  $-10$, $-6$, $-4$, $-2$, $0$, $2$, $4$, $6$, $10$\\ 
$10_{43} = \mathrm{K}(27/73):$  $-10$, $-6$, $-4$, $0$, $4$, $6$, $10$\\ 
$10_{44} = \mathrm{K}(30/79):$  $-8$, $-4$, $0$, $2$, $4$, $6$, $12$\\ 
$10_{45} = \mathrm{K}(34/89):$  $-10$, $-6$, $-4$, $-2$, $0$, $2$, $4$, $6$, $10$\\ 
$10_{46} = \mathrm{K}(1/5, 1/3, 1/2):$  $-4$, $0$, $2$, $6$, $8$, $12$, $14$, $16$\\ 
$10_{47} = \mathrm{K}(1/5, 2/3, 1/2):$  $-6$, $-2$, $0$, $4$, $6$, $10$, $12$, $14$\\ 
$10_{48} = \mathrm{K}(4/5, 1/3, 1/2):$  $-10$, $-6$, $-4$, $-2$, $0$, $2$, $6$, $8$, $10$\\ 
$10_{49} = \mathrm{K}(4/5, 2/3, 1/2):$  $-20$, $-16$, $-14$, $-12$, $-10$, $-8$, $-4$, $-2$, $0$\\ 
$10_{50} = \mathrm{K}(3/7, 1/3, 1/2):$  $-4$, $0$, $2$, $6$, $8$, $10$, $12$, $14$, $16$\\ 
$10_{51} = \mathrm{K}(3/7, 2/3, 1/2):$  $-6$, $-2$, $0$, $4$, $6$, $8$, $10$, $12$, $14$\\ 
$10_{52} = \mathrm{K}(4/7, 1/3, 1/2):$  $-10$, $-6$, $-4$, $0$, $2$, $6$, $8$, $10$\\ 
$10_{53} = \mathrm{K}(4/7, 2/3, 1/2):$  $-20$, $-16$, $-14$, $-12$, $-10$, $-8$, $-4$, $-2$, $0$\\ 
$10_{54} = \mathrm{K}(2/7, 1/3, 1/2):$  $-10$, $-6$, $-4$, $0$, $2$, $4$, $6$, $8$, $10$\\ 
$10_{55} = \mathrm{K}(2/7, 2/3, 1/2):$  $-20$, $-16$, $-14$, $-10$, $-8$, $-6$, $-4$, $-2$, $0$\\ 
$10_{56} = \mathrm{K}(5/7, 1/3, 1/2):$  $-4$, $0$, $2$, $4$, $6$, $8$, $12$, $14$, $16$\\ 
$10_{57} = \mathrm{K}(5/7, 2/3, 1/2):$  $-6$, $-2$, $0$, $2$, $4$, $6$, $10$, $12$, $14$\\ 
$10_{58} = \mathrm{K}(2/5, 2/5, 1/2 ):$  $-12$, $-8$, $-4$, $-2$, $0$, $2$, $4$, $6$, $8$\\ 
$10_{59} = \mathrm{K}(2/5, 3/5, 1/2):$  $-8$, $-4$, $0$, $2$, $4$, $6$, $8$, $10$, $12$\\ 
$10_{60}\error = \mathrm{K}(3/5, 3/5, 1/2):$  $-12$, $-8$, $-4$, $0$, $2$, $6$, $8$\\ 
$10_{61} = \mathrm{K}(1/4, 1/3, 1/3):$  $-8$, $-2$, $0$, $6$, $8$, $12$\\ 
$10_{62} = \mathrm{K}(1/4, 1/3, 2/3):$  $-6$, $0$, $2$, $8$, $10$, $14$\\ 
$10_{63} = \mathrm{K}(1/4, 2/3, 2/3):$  $-20$, $-14$, $-12$, $-8$, $-6$, $-4$, $0$\\ 
$10_{64}\errortwo = \mathrm{K}(3/4, 1/3, 1/3):$  $-8$, $-2$, $0$, $2$, $4$, $6$, $8$, $12$\\ 
$10_{65} = \mathrm{K}(3/4, 1/3, 2/3):$  $-6$, $0$, $2$, $4$, $8$, $10$, $14$\\ 
$10_{66} = \mathrm{K}(3/4, 2/3, 2/3):$  $-20$, $-14$, $-12$, $-10$, $-8$, $-6$, $-4$, $0$\\ 
$10_{67} = \mathrm{K}(2/5, 1/3, 2/3):$  $-16$, $-12$, $-10$, $-8$, $-4$, $-2$, $0$, $4$\\ 
$10_{68}\errortwo = \mathrm{K}(3/5, 1/3, 1/3):$  $-14$, $-10$, $-8$, $-6$, $-4$, $-2$, $0$, $2$, $4$, $6$\\ 
$10_{69} = \mathrm{K}(3/5, 2/3, 2/3):$  $-6$, $-2$, $0$, $2$, $4$, $6$, $8$, $10$, $12$, $14$\\ 
$10_{70} = \mathrm{K}(2/5, 1/3, 3/2):$  $-8$, $-4$, $-2$, $0$, $2$, $4$, $6$, $8$, $10$, $12$\\ 
$10_{71} = \mathrm{K}(2/5, 2/3, 3/2):$  $-10$, $-6$, $-4$, $-2$, $0$, $2$, $4$, $6$, $8$, $10$\\ 
$10_{72} = \mathrm{K}(3/5, 1/3, 3/2):$  $-4$, $0$, $2$, $4$, $6$, $8$, $10$, $12$, $14$, $16$\\ 
$10_{73} = \mathrm{K}(3/5, 2/3, 3/2):$  $-14$, $-10$, $-8$, $-6$, $-4$, $-2$, $0$, $2$, $4$, $6$\\ 
$10_{74} = \mathrm{K}(1/3, 1/3, 5/3):$  $-16$, $-10$, $-4$, $0$, $4$\\ 
$10_{75} = \mathrm{K}(2/3, 2/3, 5/3):$  $-8$, $-2$, $0$, $4$, $8$, $12$\\ 
$10_{76} = \mathrm{K}(1/3, 1/3, 5/2):$  $-4$, $0$, $2$, $6$, $8$, $12$, $14$, $16$\\ 
$10_{77} = \mathrm{K}(1/3, 2/3, 5/2):$  $-6$, $-2$, $0$, $4$, $6$, $10$, $12$, $14$\\ 
$10_{78} = \mathrm{K}(2/3, 2/3, 5/2):$  $-16$, $-12$, $-10$, $-6$, $-4$, $0$, $2$, $4$\\ 
$10_{124} = \mathrm{K}(1/5, 1/3, -1/2):$  $0$, $15$\\ 
$10_{125} = \mathrm{K}(1/5, 2/3, -1/2):$  $-10$, $0$, $4$, $32/5$\\ 
$10_{126} = \mathrm{K}(4/5, 1/3, -1/2):$  $-14$, $-8$, $-6$, $-4$, $0$, $8/3$\\ 
$10_{127} = \mathrm{K}(4/5, 2/3, -1/2):$  $-16$, $-10$, $-8$, $-6$, $-2$, $0$, $1$\\ 
$10_{128} = \mathrm{K}(3/7, 1/3, -1/2):$  $0$, $8$, $32/3$, $16$\\ 
$10_{129}\error = \mathrm{K}(3/7, 2/3, -1/2):$  $-10$, $-2$, $0$, $1$, $4$, $20/3$\\ 
$10_{130} = \mathrm{K}(4/7, 1/3, -1/2):$  $-14$, $-8$, $-4$, $0$, $8/3$\\ 
$10_{131} = \mathrm{K}(4/7, 2/3, -1/2):$  $-16$, $-10$, $-6$, $-2$, $0$, $1$\\ 
$10_{132}\error = \mathrm{K}(2/7, 1/3, -1/2):$  $-14$, $-2$, $0$, $3/2$\\ 
$10_{133} = \mathrm{K}(2/7, 2/3, -1/2):$  $-16$, $-6$, $-10/3$, $-2$, $0$, $1/2$\\ 
$10_{134} = \mathrm{K}(5/7, 1/3, -1/2):$  $0$, $4$, $6$, $8$, $10$, $14$, $50/3$\\ 
$10_{135} = \mathrm{K}(5/7, 2/3, -1/2):$  $-10$, $-6$, $-4$, $-2$, $0$, $4$, $6$, $7$\\ 
$10_{136}\error = \mathrm{K}(2/5, 2/5, -1/2):$  $-8$, $0$, $1$, $14/3$, $8$\\ 
$10_{137} = \mathrm{K}(2/5, 3/5, -1/2):$  $-12$, $-8$, $-4$, $0$, $1$, $14/3$\\ 
$10_{138}\error = \mathrm{K}(3/5, 3/5, -1/2):$  $-8$, $-4$, $-2$, $0$, $4$, $8$, $9$\\ 
$10_{139}\errortwo = \mathrm{K}(1/4, 1/3, -2/3):$  $0$, $12$, $13$, $18$, $20$\\ 
$10_{140} = \mathrm{K}(1/4, 1/3, -1/3):$  $-14$, $0$, $8/5$\\ 
$10_{141}\errortwo = \mathrm{K}(1/4, 2/3, -1/3):$  $-12$, $-4$, $-2$, $0$, $2$, $9/2$\\ 
$10_{142} = \mathrm{K}(3/4, 1/3, -2/3):$  $0$, $8$, $12$, $16$\\ 
$10_{143}\errortwo = \mathrm{K}(3/4, 1/3, -1/3):$  $-14$, $-8$, $-6$, $-2$, $0$, $8/3$\\ 
$10_{144} = \mathrm{K}(3/4, 2/3, -1/3):$  $-12$, $-8$, $-4$, $-2$, $0$, $2$, $5$\\ 
$10_{145}\errortwo = \mathrm{K}(2/5, 1/3, -2/3):$  $-18$, $-6$, $-4$, $0$, $2$\\ 
$10_{146}\errortwo = \mathrm{K}(2/5, 2/3, -1/3):$  $-10$, $-4$, $-2$, $0$, $2$, $3$, $4$, $20/3$\\ 
$10_{147}\errortwo = \mathrm{K}(3/5, 1/3, -1/3):$  $-8$, $-4$, $-2$, $0$, $4$, $6$, $26/3$\\ 

\end{semitable}

The program computes the Euler characteristic and number of boundary
components of the ``simplest'' surface for each boundary slope.  In
the above knots, the only genus zero incompressible surfaces are the
annuli in the torus knots $3_1, 5_1, 7_1, 8_{19}, 9_1$, and $10_{124}$.
The non-two-bridge Montesinos knots above with boundary slopes realized
by genus one incompressible surfaces are:

%
%

\begin{semitable}
$8_5$:  $12$ \\
$8_{20}$:     $0$ (non-Seifert surface)\\
$9_{35}$:    $0$ (Seifert surface)\\
$9_{42}$: $6$ \\
$9_{46}$: $2$, $0$ (Seifert surface)\\
$10_{46}$:        $16$\\
$10_{61}$:        $12$\\
$10_{125}$: $4$ \\
$10_{126}$:       $-4$ \\
$10_{132}$\error: $-2$\\
$10_{139}$: $12$, $13$\\
$10_{140}$: $0$ (non-Seifert surface)\\
$10_{142}$: $12$ \\
$10_{145}$: $-6$, $-4$\\
\end{semitable}

Some more complicated examples:

\begin{semitable}
$\mathrm{K}(2/5, 3/7, -1/3, -5/8)\error:$  $-14$, $-10$, $-8$, $-6$, $-16/3$, $-4$, $-23/6$, $-2$, $-8/23$, $0$, $2$, $4$, $6$, $8$, $10$, $12$, $14$, $16$, $20$, $\infty$\\ 
 \ \\ 
$\mathrm{K}(2/3, 1/3, -3/5, -3/4, 3/7):$  $-18$, $-14$, $-12$, $-10$, $-8$, $-6$, $-4$, $-2$, $-4/3$, $0$, $1/2$, $80/51$, $2$, $24/7$, $38/11$, $4$, $16/3$, $6$, $8$, $10$, $12$, $14$, $16$, $20$, $\infty$\\ 
 \ \\ 
$\mathrm{K}(1/3, 1/3, -1/3, -2/5, 1/5, -3/4, 2/3):$  $-16$, $-12$, $-10$, $-8$, $-6$, $-4$, $-2$, $0$, $2$, $4$, $6$, $8$, $10$, $58/5$, $12$, $40/3$, $122/9$, $124/9$, $14$, $102/7$, $190/13$, $76/5$, $168/11$, $142/9$, $16$, $434/27$, $146/9$, $18$, $20$, $22$, $24$, $26$, $28$, $30$, $34$, $\infty$\\ 
 \ \\ 
$\mathrm{K}(-15/32, 3/11, 7/41)\errortwoalone:$  $-34$, $-30$, $-28$, $-26$, $-24$, $-22$, $-20$, $-18$, $-16$, $-14$, $-12$, $-10$, $-8$, $-6$, $-74/15$, $-4$, $-26/7$, $-2$, $-13/8$, $-16/17$, $-10/11$, $0$, $2$, $44/19$, $40/13$, $34/11$, $74/21$, $4$, $6$, $86/11$, $8$, $148/17$, $10$, $83/8$, $152/13$, $12$, $127/10$, $216/17$, $14$, $272/19$, $16$, $167/10$, $18$, $20$, $22$, $24$\\ 
 \ \\ 
$\mathrm{K}(11/53, 17/43, -13/21):$  $-36$, $-32$, $-28$, $-24$, $-47/2$, $-22$, $-62/3$, $-20$, $-39/2$, $-18$, $-390/23$, $-50/3$, $-16$, $-44/3$, $-594/41$, $-72/5$, $-14$, $-68/5$, $-12$, $-10$, $-48/5$, $-8$, $-13/2$, $-6$, $-28/5$, $-26/5$, $-14/3$, $-22/5$, $-4$, $-5/2$, $-2$, $-6/5$, $-2/3$, $-2/5$, $0$, $2$, $3$, $7/2$, $4$, $24/5$, $6$, $15/2$, $8$, $44/5$, $10$, $23/2$, $12$, $38/3$, $64/5$, $14$, $16$, $50/3$, $18$, $20$, $22$, $24$\\ 
 \ \\ 
$\mathrm{K}(1/3, 3/5, -3/4, -2/7, 3/11, -5/13)\error:$  $-14$, $-10$, $-8$, $-6$, $-4$, $-2$, $0$, $2$, $4$, $664/117$, $6$, $20/3$, $38/5$, $8$, $62/7$, $9$, $19/2$, $776/81$, $48/5$, $260/27$, $10$, $32/3$, $98/9$, $11$, $58/5$, $82/7$, $12$, $110/9$, $112/9$, $90/7$, $13$, $27/2$, $122/9$, $68/5$, $96/7$, $14$, $72/5$, $44/3$, $134/9$, $15$, $46/3$, $108/7$, $78/5$, $110/7$, $16$, $146/9$, $148/9$, $118/7$, $17$, $52/3$, $35/2$, $88/5$, $230/13$, $124/7$, $18$, $92/5$, $170/9$, $96/5$, $58/3$, $136/7$, $138/7$, $20$, $106/5$, $64/3$, $43/2$, $282/13$, $152/7$, $22$, $116/5$, $24$, $276/11$, $26$, $28$, $30$, $32$, $34$, $36$, $38$, $42$, $\infty$\\ 

\end{semitable}

\section{Preparing the table}\label{method}

In this section I briefly describe the precautions taken to insure
that the table in Section~\ref{table} is correct.  I first noticed
errors in the table in \cite{HatcherOertel} when comparing it with the
output of a program that computes the {\em normal} boundary slopes of
a knot (that is, the boundary slopes of all surfaces which are normal
with respect to a choice of triangulation of the exterior of the
knot).  The data in Section~\ref{table} is consistent with the normal
boundary slope data.  Following the algorithm given in
\cite{HatcherOertel}, but without reference to the program used in
computing the table there, I wrote a completely new program to compute
boundary slopes.  I then compared its output with the output of
Hatcher and Oertel's program.  When the output differed, I debugged
both programs until I found the source of the problem, and then fixed
the appropriate program.  Eventually, after the new program and the
fixed version of Hatcher and Oertel's program had agreed for thousands
of trial Montesinos knots, I declared victory and went home.

Unfortunately, there remained several cases where I implemented the
algorithm incorrectly and Hatcher and Oertel's program also made a
very similar mistake.  Comparing the 1999 version of this note with
Marc Culler's computations of $A$-polynomials \cite{Culler}, Thomas
Mattman found some additional errors, leading to this 2010 revision.
For the table of knots with at most 10 crossings, the 1999 version of
this note is strictly better than \cite{HatcherOertel}, i.e.~some
errors were corrected but no new ones were introduced.  Hopefully, the
tables are now completely correct, but it is possible that errors
remain.

One source of problems in the table in \cite{HatcherOertel} is a
slight error in the main body of \cite{HatcherOertel}.  The remark
immediately preceding Proposition~2.7 gives conditions that are
claimed to be equivalent to certain hypotheses of Propositions~2.6 and
2.7.  It was this reformulation of these propositions which was used
in Hatcher and Oertel's program.  However, the conditions given in the
remark are not in fact equivalent to those given in the propositions.  My
program used the original formulation.

\section{Getting the program}\label{gettingprogram}

The programs I used in preparing this note are available at
\url{http://dunfield.info/montesinos}.  They are written the
programming language Python, and you will need a Python interpreter to
run them.  These interpreters are available, for free and for almost
all platforms, from \url{http://python.org}.

\iftopology
\else
\vspace{-0.1in}
\fi

\end{document}